# On the High Water Mark Convergents of Champernowne's Constant in Base Ten

*John K. Sikora*

## Abstract

In this paper we show that numerous patterns exist in the properties of the convergents formed by truncating the Continued Fraction Expansion (CFE) of Champernowne's Constant in base 10 ($C_{10}$) immediately before the High Water Marks (HWMs). From these patterns, we have formulated conjectures that may be used to predict the first position in $C_{10}$ that is incorrect as calculated by the convergent, the error of the convergent, the denominator of the convergent, and the number of digits of $C_{10}$ that are required to give the correct numerator of the convergent. The numerator and denominator then provide a very efficient method for calculating the CFE coefficients. In addition, we have formulated a conjecture that predicts the exact length in decimal digits of each of the HWM coefficients of $C_{10}$ for HWM numbers greater than 3. Furthermore, it is conjectured that the lengths of certain other large coefficients that appear in the CFE of $C_{10}$ are related to the lengths of the HWMs and a formula is given for this relationship. We designate these other large coefficients $2^{nd}$ Generation HWMs, or Child HMWs. It is shown that $3^{rd}$, $4^{th}$, etc., Generation HWMs also exist. A conjecture regarding the structure and exact length of the denominators of the $2^{nd}$ Generation HWM Convergents is presented. Finally, a table is presented that shows the Coefficient Numbers and the lengths of the $3^{rd}$ and $4^{th}$ Generation coefficients up to but not including HWM #12.

## Introduction

The Constant $C_{10}$ was formulated in 1933 by British mathematician and economist D. G. Champernowne as an example of a normal number [1]. It is formed by concatenating the positive integers to the right of a decimal point, i.e., 0.123456789101112…, without end. In 1937, Kurt Mahler proved that $C_{10}$ is transcendental [2].

The method of the construction of $C_{10}$, and the fact that it is normal may suggest that the Constant is somewhat mundane. However, there are two aspects of the Constant that suggest otherwise. First, any pattern of digits, no matter the content or length, will eventually appear in $C_{10}$ [3]. In fact, the location of the pattern (although not necessarily the first one) may be calculated [4][5]. Second, and the focus of this paper, is the nature of the CFE coefficients of $C_{10}$. The CFE consists mostly of coefficients with a reasonably small number of digits, interspersed with coefficients with a very large number of digits. A coefficient that has more digits than any previous coefficient is called a HWM. As one progresses along the CFE, the pattern of relatively small coefficients sprinkled with very large coefficients continues between the HWMs, even though the other large coefficients are not themselves HWMs.

At present, the HWMs are known up to and including HWM #11, which is coefficient number 13,522 (per OEIS A143533) [6] and the coefficient has 651,111,094 digits. The number of digits in each of the known HWMs of $C_{10}$ may be found in sequence OEIS A143534 [7], as well as in Table 1. The present work has determined that the coefficient number of HWM #12 is 34,062 per the numbering system given in [6] and in this paper. This result is believed to be previously unknown.

It is well known that the HWMs in the CFE of $C_{10}$ increase without bound, and that calculating the coefficients is difficult. The purpose of this paper is to present the myriad interesting patterns in and the properties of the Convergents Before the HWMs of $C_{10}$, formulate conjectures that allow the very efficient calculation of the CFE coefficients of $C_{10}$, and to encourage research, both theoretical and applied, into $C_X$, where $C_X$ denotes the Champernowne Constant in base X, including X = 10.

# Outline

This paper consists of a number of sections. First, the terminology that is used in the paper is given, and then the methods of analysis that were used are presented. Section 1 presents the conjectures that resulted from an analysis of the characteristics of the HWM Convergents, with the HWM Convergents as defined below. These conjectures are used to point out the interesting and useful patterns in the characteristics of the HWM Convergents. Section 1 includes a table of these patterns and characteristics.

Section 2 shows how the conjectures may be used to efficiently calculate the $C_{10}$ CFE coefficients up to each HWM, for HWM Numbers greater than 4. The special case of HWM #4 is also stated.

Section 3 examines the lengths of the 2nd Generation HWM coefficients, the characteristics of the Convergents of the $2^{nd}$ Generation HWMs, and their possible use in the calculation of the $C_{10}$ CFE coefficients. Currently, this is not possible, as a part of the pattern required for the calculation has not yet been decoded. This section tabulates the known patterns and explains what parts of the pattern are missing. The reasons for possibly using the $2^{nd}$ Generation HWMs instead of the HWMs ($1^{st}$ Generation) for the $C_{10}$ CFE coefficient calculation are discussed. In addition, the existence of $3^{rd}$, $4^{th}$, etc., Generation HWMs is demonstrated.

Section 4 explains the progression of the computational methods, the limitations we faced, and the workarounds. Access to the programs for other researchers is also discussed.

Section 5 discusses ways that the present analysis may be improved and extended. It also suggests areas for future research. Finally, the Conclusions are stated and a Summary of Findings is presented.

# Terminology

*$C_{10}$ Position, or Position* - For an integer, it is the position of the first digit of the integer in the consecutive integer sequence used to generate $C_{10}$, starting at Position 0 for the 0 to the left of the decimal point. The decimal point is not counted. Therefore, the integer 1 is in Position 1; the integer 12 is in Position 14, etc. For a specific digit of an integer, it is the position as just described, except for that specific digit. Thus, the 2 in the integer 12 is in Position 15.

*CFE Coefficient Number, or Coefficient Number* - The number of the CFE coefficient, starting at Coefficient Number 0 (for the 0). For the HWMs, this corresponds to the sequence OEIS A143533, which is OEIS sequence A038705(n) – 1 [8].

*HWM Convergent, or The Convergent Before HWM #N, or Convergent* - The convergent formed from the $C_{10}$ CFE coefficients terminated immediately before the HWM #N coefficient. This also applies to the 2nd Generation HWM Convergents, meaning that the CFE coefficients are terminated immediately before the 2nd Generation HWM term.

*Efficiency of Calculation of the CFE Coefficients* - The efficiency specifically refers to the number of digits of $C_{10}$ that are required to give a certain number of CFE coefficient digits, regardless of the method used.

*HWM Number, or HWM #* - The Number of the HWM, with the first HWM corresponding to HWM #1. Thus, the Coefficient Number of HWM #N corresponds to the Nth term in sequence OEIS A143533, starting at N = 1.

*Number of Correct Digits (NCD)* - The number of correct digits is the number of correct consecutive digits in $C_{10}$, starting at, and counting the 0 to the left of the decimal point, but not counting the decimal point itself. Therefore, if the Position of the last correct digit, $P_C$, is known, the number of correct digits is $P_C + 1$ since the 0 is counted as the first correct digit, but it is in Position 0. Similarly, if the Position of the first incorrect digit, $P_I$, is known, the number of correct digits is $P_I$.

*TCFE* – Truncated Continued Fraction Expansion.

*MPIR* - Multiple Precision Integers and Rationals Library for C

*GMP* – GNU Multiple Precision Arithmetic Library for C

*$C_{10}$ Position of $10^M$* – The formula for the $C_{10}$ Position of the integer $10^M$, $M \geq 0$, per the above definition is:

$$C_{10}\ Position\ of\ 10^M = 1 + \left( \sum_{m=1}^{M} 9 \cdot m \cdot 10^{(m-1)}, m \geq 1 \right)$$

Thus, the $C_{10}$ Position of $10^0$ is 1; the Position of $10^1$ is 10, and the Position of $10^2$ is 190, etc.

# Methods Used

The unusual properties of the CFE coefficients, particularly those of the HWM Coefficients of $C_{10}$ piqued our interest. Using the coefficient values available online up to Coefficient Number 161 [9], which is the coefficient immediately before HWM #7, a Ruby program was written to calculate the HWM Convergents Before HWM #5 through HWM #7. In each case, the program also determined the first failed digit and calculated the error (the difference from $C_{10}$). This was relatively easy as Ruby has arbitrary precision integer and floating point arithmetic functions built into the standard library.

We noticed patterns that allowed us to formulate conjectures about the succeeding HWMs. Assuming that the conjectures held true, we proceeded to calculate our own coefficients for the succeeding HWMs. The subsequent results verified our assumption. Due to speed considerations, C programs using MPIR and then GMP were written to calculate the CFE coefficients up to Coefficient Number 34061, the coefficient immediately before HWM #12. However, the error was verified only up to the Convergent Before HWM #11 (Coefficient Number 13521).

Various Ruby and C programs were written to check the results between the Ruby and the C output files, the number of digits in (the lengths of) the HWMs, the lengths of the 2nd Generation HWMs, the lengths of the other coefficients, and to calculate the various results for the Convergents Before the 2nd Generation HWMs.

Section 4 provides more detail on the methods used, especially the computational methods.

# Section 1      Patterns, Predictions, and Conjectures on the Convergents of $C_{10}$ Truncated Immediately Before the HWMs

In this section, the conjectures described earlier are stated and explained. The next section shows how these conjectures may be used to efficiently calculate the CFE coefficients up to each of the HWMs, assuming that the conjectures hold true.

It should be noted that most of the formulae in the conjectures are given in terms of the Number of Correct Digits (NCD) as calculated by the Convergents. Alternate forms of the formulae in terms of the $C_{10}$ Position of $10^n$, where n is a function of the HWM Number may be derived, but are not given in this paper.

In order to provide a preview, before the conjectures are stated we will present the following table so that the patterns in the characteristics of the Convergents Before the HWMs are more apparent:

# Table 1 - $C_{10}$ HWM CFE Characteristics Summary

| Convergent (CFE Truncated) *Before* HWM Number | Coefficient Number of the HWM | Integer in $C_{10}$ That Fails; Integer That it Fails As in the Convergent Calculation | Number of Correct Digits Calculated By the Convergent | Convergent Error | Convergent Denominator | Number of Digits in the Next Term (the HWM) |
|---|---|---|---|---|---|---|
| 4 | 4 | 8; 9 | 8 | 1.0E-9 | 81 (divided by 2 = 4.05E+1, see text) | 6 |
| 5 | 18 | 98; 99 | 187 | 9.1E-190 | 4.9005E+11 | 166 |
| 6 | 40 | 998; 999 | 2886 | 9.01E-2890 | 4.990005E+192 | 2504 |
| 7 | 162 | 9998; 9999 | 38,885 | 9.001E-38,890 | 4.99900005E+2893 | 33,102 |
| 8 | 526 | 99,998; 99,999 | 488,884 | 9.0001E-488,890 | 4.9999000005E+38,894 | 411,100 |
| 9 | 1708 | 999,998; 999,999 | 5,888,883 | 9.00001E-5,888,890 | 4.999990000005E+488,895 | 4,911,098 |
| 10 | 4838 | 9,999,998; 9,999,999 | 68,888,882 | 9.000001E-68,888,890 | 4.99999900000005E+5,888,896 | 57,111,096 |
| 11 | 13,522 | 99,999,998; 99,999,999 | 788,888,881 | 9.0000001E-788,888,890 | 4.9999999000000005E+68,888,897 | 651,111,094 |
| 12 | 34,062 | 999,999,998; 999,999,999 | 8,888,888,880 | 9.00000001E-8,888,888,890 | 4.999999990000000005E+788,888,898 | 7,311,111,092 |
| 13 | ? | 9,999,999,998; 9,999,999,999 | 98,888,888,879 | 9.000000001E-98,888,888,890 | 4.99999999900000000005E+8,888,888,899 | 81,111,111,090 |
| 14 | ? | 99,999,999,998; 99,999,999,999 | 1,088,888,888,878 | 9.0000000001E-1,088,888,888,890 | 4.9999999999000000000005E+98,888,888,900 | 891,111,111,088 |
| | Conjecture 3 (Even Numbered) | Conjecture 1 | Conjecture 1, Conjecture 2 | Conjecture 1, Conjecture 4 | Conjecture 6 | Conjecture 5 |

The lightly shaded cells have values that have been confirmed by our calculations and in some cases the values have also been confirmed by others [15]. We believe that the value (34062) in the medium shaded cell, calculated using the method presented in this paper, was previously unknown. To the best of our knowledge, the values in the darkest shaded cells are predictions, and they have not been confirmed.

The table makes it easier to recognize most, if not all, of the patterns. However, the conjectures that follow give the formulae behind the patterns. These formulae allow for the verification of the characteristics of the known HWMs, and provide the ability to predict the characteristics of, and calculate the CFE coefficients up to higher numbered HWMs.

The conjectures are as follows:

## Conjecture 1

Let the Convergent Before HWM #N be the convergent formed by the $C_{10}$ TCFE as truncated immediately before HWM #N, N ≥ 4. Then the convergent calculates the value of $C_{10}$ correct to the last 8 before the number $10^{(N-3)}$ in the integer sequence used to form $C_{10}$. Furthermore, the 8 is calculated as a 9, followed by [2·(N - 4) + 1] zeroes, followed by a 1, followed by N − 4 zeroes, then followed by a 2.

Thus, the Convergent Before HWM #4 calculates $C_{10}$ as 0.12345679012…, the Convergent Before HWM #5 calculates $C_{10}$ as 0.12345…969799000102…, etc. It should be noted that HWM #4 is the first HWM to reach an integer of the form, $10^X$, X ≥ 1, in the $C_{10}$ integer sequence. Thus, the integer $10^1 = 10$ is the first integer with a preceding 8 that is available for the convergent to fail to calculate correctly.

We have not seen this Conjecture explicitly stated. However, later study of an available Mathematica notebook [11] shows that portions of Conjecture 1 have been previously observed.

## Conjecture 2

The Number of Correct Digits including the 0 to the left of the decimal point calculated by the Convergent Before HWM #N, N ≥ 3, is:

$$Number\ of\ Correct\ Digits = NCD(N) = C_{10}\ Position\ of\ 10^{(N-3)} - N + 2$$

Thus, NCD(3) = 0, NCD(4) = 8, NCD(5) = 187, etc. Conjecture 2 follows from Conjecture 1. The Number of Correct Digits (calculated) of the Convergent Before HWM #N is explicitly given because the value appears in a number of the following conjectures and formulae. It should be noted that according to the formulae, NCD(3) = 0. This is convenient in that it allows for the formulae that follow to be true for all HWM #N, N ≥ 4, as there are some NCD(N – 1) terms. The Convergent Before HWM #4 is the first one that calculates correct digits up to the 8 before the first positive integer exponent of 10 ($10^1$) in the sequence used to form $C_{10}$. Therefore, in this context it makes sense to have a formula where NCD(3) = 0.

## Conjecture 3

The HWM Coefficient Numbers starting at HWM #4 must be even.

Note that the Coefficient Numbers are as defined above and are as given in OEIS A143533. Conjecture 3 follows from the statement in Conjecture 1 that the last 8 before the power of 10 boundary is calculated as a 9 by the Convergent. With the Coefficient Numbers starting at 0 instead of 1, the odd numbered coefficients will calculate numbers slightly larger than $C_{10}$. Therefore, the CFE must terminate on an odd numbered coefficient. Since the CFE is terminated immediately before the HWM, the HWM coefficient numbers after HWM #3 must be even, per the stated numbering scheme.

This conjecture is very useful at the end of the calculation of the CFE coefficients as the algorithm cannot distinguish the degenerate case where the last coefficient is 1 from the case where the last coefficient is not 1. In the degenerate case where the last coefficient equals 1, the algorithm calculates the last term as $Y + 1$ instead of $Y$, and the calculation ends on an even numbered coefficient. However, if the calculation ends on an even numbered coefficient, it is known from Conjecture 3 that this is the degenerate case, the even numbered coefficient is $Y$, and the last coefficient (the one before the HWM) is 1. Otherwise, the HWM would occur at an odd numbered coefficient. Section 2 describes the CFE calculation process in more detail.

## Conjecture 4

The error of the Convergent Before HWM #4 is approximately 1.0E-9.  Let NCD(N) be the Number of Correct Digits calculated by the Convergent Before HWM #N, N ≥ 5. Then the error of the Convergent  Before HWM #N is approximately:

$$9.\underbrace{0001}\text{E-exp}$$

Number of
Zeroes

The Number of Zeroes is $N - 5$, and exp (the exponent) is equal to $NCD(N) + N - 2$ (which is equal to the $C_{10}$ Position of $10^{(N-3)}$).

Conjecture 4 follows as a consequence of Conjecture 1. For example, the error of the Convergent Before HWM #5 is approximately 9.1E-190, and the error of the Convergent Before HWM #6 is approximately 9.01E-2890, etc.

The magnitude of the error exponent was instrumental in the course of this study. Specifically, it was the fact that an 11 digit numerator divided by a 12 digit denominator (from the 18 CFE terms before HWM #5) results in an error from $C_{10}$ within a few orders of magnitude of a googol[th] of a googol[th] of a unit that led us to further investigation of the Convergents' error and characteristics.

## Conjecture 5

Let NCD(N) be the Number of Correct Digits calculated by the Convergent Before HWM #N. Then the length in decimal digits of HWM #N, N ≥ 4, is:

$$Length\ of\ HWM\ \#N = NCD(N) - 2 \cdot NCD(N-1) - 3 \cdot (N-2) + 4$$

For example, the length of HWM #4 is 6 digits, the length of HWM #5 is 166 digits, and the length of HWM #6 is 2504 digits. The largest known HWM is HWM #11 which has 651,111,094 digits. Note that NCD(3) = 0. This comes from the equation for NCD(N), and it also stands to reason from the fact that the Convergent Before HWM #4 is the first one that calculates correct digits to the 8 before the first positive integer exponent of 10 $(10^1)$ in the sequence used to form $C_{10}$, which is the premise of the conjectures.

It is predicted that the exact lengths of HWM #12, HWM #13, and HWM #14 of $C_{10}$ are 7,311,111,092; 81,111,111,090; and 891,111,111,088 digits, respectively, following the formula given in Conjecture 5.

## Conjecture 6

Let NCD(N) be the Number of Correct Digits calculated by the Convergent Before HWM #N. Then the denominator of the Convergent Before HWM #N, for N ≥ 5, is:

4.99900005E+exp

Number of          Number of
Nines              Zeroes

The Number of Nines is N − 4, the Number of Zeroes is N - 3, and exp (the exponent) is:

$$\exp = NCD(N-1) + 2 \cdot (N-2) - 3$$

For example, the denominator for the Convergent Before HWM #5 is 4.9005E+11, and the denominator for the Convergent Before HWM #6 is 4.990005E+192.

Note that the denominator of the Convergent Before HWM #4 is 81, the numerator is 10, and that:

$$\frac{10}{81} = \frac{10/2}{81/2} = \frac{5}{40.5}$$

When the numerator and denominator are divided by 2, the denominator becomes 40.5, or 4.05E+1. Since NCD(3) = 0, per the conjecture, the exponent is +1, and there are 4 - 4 = 0 '9's, and 4 – 3 = 1 '0'. Thus, 4.05E+1 satisfies the formula given in Conjecture 6 for HWM #4 if a non-integer simplification is allowed.

Conjecture 6 is very useful, as the knowledge of the exact value of the denominator allows for the calculation of the numerator with a sufficient, but minimal number of digits of $C_{10}$. The CFE coefficients of $C_{10}$ may then be calculated from the numerator and denominator. Section 2 describes this process in more detail.

The output files verifying the calculations up to the Convergent Before HWM #11 may be found in [18], along with the program source code. The file containing the coefficients calculated from the Convergent Before HWM #12 is about 789 MB (unzipped), and may be found in [19]. Section 4 provides more details.

# Section 2      Efficient Calculation of the CFE Coefficients of $C_{10}$ Given the Preceding Conjectures on the HWM Convergents

The $C_{10}$ CFE Coefficients may be calculated very efficiently from the Convergents terminated immediately before a HWM. Since the denominator of the convergent is known, the numerator of the convergent may be calculated by multiplying the denominator by a sufficient number of digits of $C_{10}$.

By using the word 'efficient', we specifically mean that the number of digits of $C_{10}$ that is required to calculate the resulting number of coefficients is kept as low as possible. This process for calculating the TCFE coefficients is also efficient in the sense that once the numerator is found the calculations are performed on integers (or rational numbers) and not on floating point numbers. This results in high efficiency from a calculation time aspect.

The high efficiency of the process is a result of the following conjecture:

## Conjecture 7

Given:

1) The denominator, $D_N$, of a Convergent truncated immediately before HWM #N, N ≥ 5, as calculated from Conjecture 6

2) The $C_{10}$ Position, $P_N$, of the integer $10^{(N-4)}$

3) $C_{10}[0..P_N]$ denotes $C_{10}$ truncated to Position $P_N$

4) Ceil(F) denotes the integer immediately above non-integer, F

Then the numerator, $N_N$, of the Convergent is calculated as:

$$N_N = Ceil(D_N \cdot C_{10}[0..P_N])$$

For example, for the Convergent Before HWM #5:

$D_5 = 4.9005\text{E}+11 = 490{,}050{,}000{,}000$

$P_5$ = the $C_{10}$ Position of the integer $10^{(5-4)} = 10^1 = 10$, the Position of which itself is 10, so $P_5 = 10$

$C_{10}[0..P_5] = C_{10}[0..10] = 0.1234567891$

and

$$N_5 = Ceil(490{,}050{,}000{,}000 \cdot 0.1234567891) = Ceil(60499999498.455)$$
$$= 60{,}499{,}999{,}499$$

Note that for HWM #4, the numerator is 81. If it is divided by two the result is 40.5 (which fits the pattern for the denominator given in Conjecture 6). The $C_{10}$ Position of $10^0$, $P_4$, is 1 so $C_{10}[0..P_4] = 0.1$. The numerator may then be calculated as Ceil(40.5 · 0.1) = 5. If both the numerator and denominator are then multiplied by two to get the original denominator, the correct Convergent, $^{10}/_{81}$ is calculated.

Conjecture 6 defines the denominator, Conjecture 7 allows for the calculation of the numerator, and Conjecture 1, Conjecture 2, and Conjecture 4 describe the result of the numerator divided by denominator of the Convergent to give the approximation to $C_{10}$. If the previous conjectures are taken to be true, then Conjecture 7 can be shown to be true. However, this is beyond the scope of this paper.

Once the numerator is calculated, the CFE coefficients may be calculated using the well-known algorithm [12] from the numerator and the denominator. When the calculation has reached the final coefficient, Conjecture 3 is used to determine if the last coefficient is 1. If it is, then the second to last coefficient is equal to the calculated value minus 1.

## Efficiency of the Method

In order to show the efficiency of the method, we first compare an example calculation to perhaps the least efficient method of calculation of the CFE for a non-rational number such as $C_{10}$. In this method, a truncated version of $C_{10}$ is used. The first coefficient is the integer value of $C_{10}$ to the left of the decimal point. The reciprocal of the remainder is then taken. The integer value of the result is the next coefficient and the process continues by taking the reciprocal of the remainder, and repeating.

Using the process described in the paragraph above, if we take the sequence of $C_{10}[0..P_5] = 0.1234567891$ as determined from Conjecture 7, the CFE coefficients of $C_{10}$ are calculated as [0; 8, 9, 1, 148921, …]. The last term is not correct, as it should be 149083. Thus, this method calculates only 6 correct digits (including the 0).

By contrast, if we use $C_{10}[0..P_5] = 0.1234567891$ as determined by Conjecture 7 and the denominator of 490,050,000,000 from Conjecture 6, we calculate a numerator of 60,499,999,499. If we then calculate the CFE coefficients from the numerator and denominator (using a process similar to that described above, except with integers), we get the coefficients correctly as [0; 8, 9, 1, 149083, 1, 1, 1, 4, 1, 1, 1, 3, 4, 1, 1, 1, 15]. Thus, 24 digits are correctly calculated using 11 digits of $C_{10}$.

The following table shows the confirmed efficiency statistics of our calculations of the CFE coefficients from the HWM Convergents:

**Table 2 - CFE Coefficient Calculation Efficiency Statistics**

| Convergent Before HWM Number | Total Number of CFE Digits Calculated | Number of Digits of $C_{10}$ Required for the Calculation | $C_{10}[0..P_N]$ Used in the Calculation |
|---|---|---|---|
| 4 | 4 | 2 | 0.1 |
| 5 | 24 | 11 | 0.1234567891 |
| 6 | 217 | 191 | 0.12345…98991 |
| 7 | 2995 | 2891 | 0.12345…9989991 |
| 8 | 39,231 | 38,891 | 0.12345…999899991 |
| 9 | 489,981 | 488,891 | 0.12345…99998999991 |
| 10 | 5,891,976 | 5,888,891 | 0.12345…9999989999991 |
| 11 | 68,897,496 | 68,888,891 | 0.12345…999999899999991 |
| 12 | 788,910,618 | 788,888,891 | 0.12345…99999998999999991 |
| | | **Part of Conjecture 7** | |

The Convergent numerator and denominator are guaranteed to be in lowest terms [13] so that common divisors are not present. Thus, the calculation involves the smallest numerator and denominator possible, maintaining high efficiency.

The disadvantage of the method is that each successive HWM requires an increase of approximately 12 times in memory and an increase of approximately 24 times in computing time.

**Other Efficient Methods of Calculation of the CFE Coefficients**

There are other efficient methods for the calculation of the CFE coefficients of $C_{10}$ [10][11]. However, we believe that the method presented in this paper is the most efficient that has presently been found. Using the method given in [10], in order to calculate through Coefficient Number 38 (215 total coefficient digits), 383 decimal digits of $C_{10}$ are required. Using the method described in this paper, only 191 decimal digits of $C_{10}$ are required to calculate the CFE coefficients through Coefficient Number 39 (217 coefficient digits).

In [11], it is stated that in order to get the first 160 coefficients using the method given in the reference (2993 coefficient digits), 5787 decimal digits of $C_{10}$ are required. It is also noted that the numerator and the denominator as calculated by this method require reduction to lowest terms. Using the method described in this paper, calculating the first 162 coefficients (2995 coefficient digits), requires only 2891 decimal digits of $C_{10}$, and the numerator and denominator are in lowest terms, thus reducing the calculation time.

The drawback in using the method presented in this paper is that the coefficients may be calculated in this manner only at each successive HWM.

# Section 3 $2^{nd}$, $3^{rd}$, etc., Generation HWMs (Child HWMs)

It is well known that in addition to the HWM $C_{10}$ CFE coefficients, there are other coefficients that have a large number of digits that are not themselves HWMs. For example, coefficient 101 (starting from 0) is 140 digits long, but it is not a HWM since coefficient number 40 is 2504 digits long and coefficient number 18 is 166 digits long. Coefficient 101 is an example of what we are designating a $2^{nd}$ Generation HWM, or a Child HWM.

## $2^{nd}$ Generation HWM Lengths

A pattern has been observed in the number of digits of the $2^{nd}$ Generation HWMs of $C_{10}$, leading to the following conjecture:

## Conjecture 8

Given HWM #(N), N ≥ 5, then between HWM #(N + 1) and HWM #(N + 2), there exists one 2nd Generation HWM, whose length in decimal digits is:

$$Length\ of\ 2^{nd}\ Generation\ HWM = [Length\ of\ HWM\ \#N] - 10 \cdot (N - 5) - 26$$

Thus, between HWM #6 and HWM #7 there exists a 2nd Generation HWM whose length is 140 digits, since HWM #5 has 166 digits. Similarly, between HWM #7 and HWM #8 there exists a 2nd Generation HWM whose length is 2468 digits as HWM #6 has 2504 digits.

The following table summarizes the results:

### Table 3 - 2nd Generation HWM Lengths

| 2nd Generation HWM Coefficient Number | Appears After HWM #(N + 1) | Length of 2nd Generation HWM | Length of HWM #N | Difference (Length of HWM is Greater) |
|---|---|---|---|---|
| 101 | 6 | 140 | 166 | 26 |
| 357 | 7 | 2468 | 2504 | 36 |
| 1221 | 8 | 33,056 | 33,102 | 46 |
| 3569 | 9 | 411,044 | 411,100 | 56 |
| 9827 | 10 | 4,911,032 | 4,911,098 | 66 |
| 25,069 | 11 | 57,111,020 | 57,111,096 | 76 |
| ? | 12 | 651,111,008 | 651,111,094 | 86 |
| ? | 13 | 7,311,110,996 | 7,311,111,092 | 96 |
| | | **Conjecture 8** | **Conjecture 5** | **Conjecture 8** |

The lightly shaded cells have values that have been confirmed by our calculations. To the best of our knowledge, the values in the darkest shaded cells are predictions, and they have not been confirmed.

## Convergents Before 2nd Generation HWMs (Convergents Terminated Immediately Before the 2nd Generation HWM) - Error

If a Convergent is taken immediately before a 2nd Generation HWM, the error is predictable, but there is no consistent pattern in the failed integer that is used to construct $C_{10}$ as there is in the HWMs.

## Conjecture 9

Let NCD(N) be the Number of Correct Digits calculated by the Convergent Before HWM #N. For a 2nd Generation HWM between HWM #N and HWM #(N + 1), N ≥ 6, the error of the Convergent Before the 2nd Generation HWM is:

$$-8.9992\text{E-exp}$$

<div align="center">Number of<br>Nines</div>

The Number of Nines is N − 5, and exp (the exponent) is:

$$\exp = -2 \cdot (error\ exponent\ of\ HWM\ \#N) - NCD(N-1) - N + 3$$

The error exponent of HWM #N is taken as a negative number. The error for the 2nd Generation HWM Convergent truncated immediately before coefficient number 101, which is between HWM #6 and HWM #7 is −8.92E-5590. The error for the 2nd Generation HWM Convergent truncated immediately before coefficient number 357, which is between HWM #7 and HWM #8 is −8.992E-74,890.

Unlike the HWMs, the integer and digit in the sequence used to form $C_{10}$ that fails is not consistent for the 2nd Generation HWMs. The digit that fails may be determined by the integer and digit number whose $C_{10}$ Position Number is (exp − 1) as given above (exp is a positive integer). For example, for the Convergent Before the 2nd Generation HWM whose coefficient number is 101, the error exponent is −5590, so exp = 5590. The integer and digit number of $C_{10}$ Position 5590 – 1 = 5589 is integer 1634, digit 4. Indeed, this Convergent calculates $C_{10}$ as 0.12345…167216731673, thus the 4th digit of integer 1674 is calculated as a 3 rather than a 4. The integer and digit number, given a $C_{10}$ Position, may be calculated using the algorithm given in [14].

Conjecture 9 indicates that the error calculated by the Convergent Before each 2nd Generation HWM is negative. Therefore, this implies that the CFE Coefficient Number (as defined earlier in this paper), is odd.

The following table summarizes the 2nd Generation HWM error results:

**Table 4 - 2nd Generation HWM Error Statistics**

| 2nd Generation HWM Appears After HWM Number | Convergent CFE Truncated before 2nd Generation HWM Coefficient Number | Convergent Error | Integer in $C_{10}$ That Fails; Integer That it Fails As in TCFE Calculation |
|---|---|---|---|
| 6 | 101 | -8.92E-5590 | 1674;1673 |
| 7 | 357 | -8.992E-74,890 | 17,199;17,198 |
| 8 | 1221 | -8.9992E-938,890 | 174,999;174,998 |
| 9 | 3569 | -8.99992E-11,288,890 | 1,771,428; 1,770,528 |
| 10 | 9827 | -8.999992E-131,888,890 | 17,874,999; 17,874,998 |
| 11 | 25,069 | -8.9999992E-1,508,888,890 | 179,999,999; 179,999,998 |
| 12 | ? | -8.99999992E-16,988,888,890 | 1,809,999,999;1,809,999,998 |
| 13 | ? | -8.999999992E-188,888,888,890 | 18,181,818,181; 18,181,818,xyz |
| | Conjecture 9 (Odd Numbered) | Conjecture 9 | [14] |

The lightly shaded cells have values that have been confirmed by our calculations. To the best of our knowledge, the values in the darkest shaded cells are predictions, and they have not been confirmed.

## Convergents Before 2nd Generation HWMs – Denominator and Potential CFE Coefficient Calculation

Since the Convergents of the CFEs truncated immediately before the HWMs yield denominators with values that may be determined via conjectures as a result of patterns that are observed, and since the denominators can be used to efficiently calculate the $C_{10}$ CFE coefficients, it is natural to ask if the Convergents truncated immediately before the 2nd Generation HWMs may similarly be used to calculate the CFE coefficients.

The advantage of using this approach would be that less memory would be needed since the calculation would extend only to the 2nd Generation HWM and not to the next HWM. Unfortunately, at present, the denominators of the Convergents Before the 2nd Generation HWMs have only partially predictable patterns.

It has been found that the denominators of this type have the following decimal form:

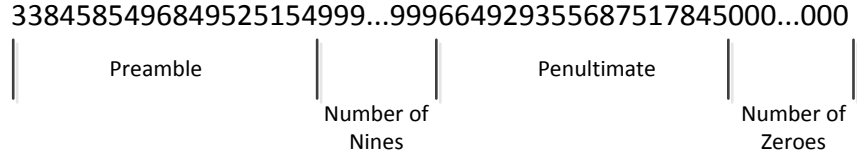

3384585496849525154999...99966492935568751784 5000...000

Preamble | Number of Nines | Penultimate | Number of Zeroes

## Conjecture 10

Let NCD(N) be the Number of Correct Digits calculated by the Convergent Before HWM #N, $N \geq 6$. The denominator $D_{N\_2}$ of the Convergent of the $2^{nd}$ Generation HWM that appears after HWM #N of $C_{10}$ consists of a Preamble $PR_N$, followed by a number of nines, $N_{N\_nines}$, followed by a Penultimate, $PE_N$, followed by a number of zeroes, $Z_N$. The length of the Preamble, $PR_N$, is given by:

$$Length \ of \ PR_N = 7 \cdot (N - 2) - 9$$

The Number of Nines, $N_{N\_nines}$, is given by:

$$Number \ of \ Nines = \ N_{N\_nines} = NCD(N) - NCD(N - 1) - 7 \cdot (N - 2) + 10$$

The length of the Penultimate, $PE_N$, is given by:

$$Length \ of \ PE_N = Length \ of \ PR_N - 1$$

The Number of Zeroes, $Z_N$, is given by:

$$Number \ of \ Zeroes = Z_N = NCD(N - 2) - 1$$

Thus, the total number of digits of $D_{N\_2}$ is given by:

$$Total \ Number \ of \ Digits \ of \ D_{N\_2}$$
$$= NCD(N) - NCD(N - 1) + NCD(N - 2) + 7 \cdot (N - 2) - 10$$

Note that there are no $2^{nd}$ Generation HWMs until after HWM #6. The total number of digits of $D_{N\_2}$ is significant because this roughly gives the number of digits of $C_{10}$ required to determine the numerator for a potential CFE calculation from the Convergents Before the $2^{nd}$ Generation HWMs. In actuality, several more digits of $C_{10}$ would be needed for this calculation. However, we have not decoded the Preamble and Penultimate content sufficiently for this purpose.

To show how striking the appearance of such a denominator can be, the following is the decimal representation of the denominator of the Convergent Before the 2[nd] Generation HWM at Coefficient Number 101 (the first and smallest one):

```
33845854968495251549999999999999999999999999999999999999999999999999999999999999
99999999999999999999999999999999999999999999999999999999999999999999999999999999
99999999999999999999999999999999999999999999999999999999999999999999999999999999
99999999999999999999999999999999999999999999999999999999999999999999999999999999
99999999999999999999999999999999999999999999999999999999999999999999999999999999
99999999999999999999999999999999999999999999999999999999999999999999999999999999
99999999999999999999999999999999999999999999999999999999999999999999999999999999
99999999999999999999999999999999999999999999999999999999999999999999999999999999
99999999999999999999999999999999999999999999999999999999999999999999999999999999
99999999999999999999999999999999999999999999999999999999999999999999999999999999
99999999999999999999999999999999999999999999999999999999999999999999999999999999
99999999999999999999999999999999999999999999999999999999999999999999999999999999
99999999999999999999999999999999999999999999999999999999999999999999999999999999
99999999999999999999999999999999999999999999999999999999999999999999999999999999
99999999999999999999999999999999999999999999999999999999999999999999999999999999
99999999999999999999999999999999999999999999999999999999999999999999999999999999
99999999999999999999999999999999999999999999999999999999999999999999999999999999
99999999999999999999999999999999999999999999999999999999999999999999999999999999
99999999999999999999999999999999999999999999999999999999999999999999999999999999
99999999999999999999999999999999999999999999999999999999999999999999999999999999
99999999999999999999999999999999999999999999999999999999999999999999999999999999
99999999999999999999999999999999999999999999999999999999999999999999999999999999
99999999999999999999999999999999999999999999999999999999999999999999999999999999
99999999999999999999999999999999999999999999999999999999999999999999999999999999
99999999999999999999999999999999999999999999999999999999999999999999999999999999
99999999999999999999999999999999999999999999999999999999999999999999999999999999
99999999999999999999999999999999999999999999999999999999999999999999999999999999
99999999999999999999999999999999999999999999999999999999999999999999999999999999
99999999999999999999999999999999999999999999999999999999999999999999999999999999
99999999999999999999999999999999999999999999999999999999999999999999999999999999
99999999999999999999999999999999999999999999999999999999999999999999999999999999
99999999999999999999999999999999999996649293556875178450000000
```

The following table summarizes the characteristics for the form of the Convergent denominators, the Convergents being truncated immediately before a 2nd Generation HWM:

**Table 5 - 2nd Generation HWM Convergent Denominator Characteristics**

| 2nd Generation HWM Appears After HWM Number | 2nd Generation HWM Coefficient Number | Preamble, $PR_N$ | Preamble Length | Number of Nines, $N_{N\_nines}$ | Penultimate, $PE_N$ | Pen-ultimate Length | Number of Zeroes, $Z_N$ | Total Length of Denominator |
|---|---|---|---|---|---|---|---|---|
| 6 | 101 | 3384585496849525154 | 19 | 2681 | 664929355687517845 | 18 | 7 | 2725 |
| 7 | 357 | 34248335240752097390773154 | 26 | 35,974 | 65785608244226627024068458 | 25 | 186 | 36,211 |
| 8 | 1221 | 342890265624321082344399571003154 | 33 | 449,967 | 657143684279973146379986357 96845 | 32 | 2885 | 452,917 |
| 9 | 3569 | 34293100142167679021625913 2673148303154 | 40 | 5,399,960 | 657072393632921948456613306 036419696845 | 39 | 38,884 | 5,438,923 |
| 10 | 9827 | 3429350754463155348149074283 8237205586421303154 | 47 | 62,999,953 | 657065264059786008368209858 8614554259258696845 | 46 | 488,883 | 63,488,929 |
| 11 | 25,069 | 342935482853228587448395154 3230851508810524691513031544 | 54 | 719,999,946 | 657064551097387983539232407 4056392318382098764869684544 | 53 | 5,888,882 | 725,888,935 |
| 12 | ? | 3429{53 more digits}3154 | 61 | 8,099,999,939 | 6570{52 more digits}6845 | 60 | 68,888,881 | 8,168,888,941 |
| | Conjecture 9 (Odd Number) | Conjecture 10 | | | | | | |

The lightly shaded cells have values that have been confirmed by our calculations. To the best of our knowledge, the values in the darkest shaded cells are predictions, and they have not been confirmed. Note that the digits that are referenced in the brackets for the Preamble and Penultimate of the 2nd Generation HWM appearing after HWM #12 are currently mostly unknown.

If formulae for the Preamble and Penultimate are discovered, then the CFE coefficients up to the 2nd Generation HWMs could be calculated by finding the numerators as discussed above. However, as can be seen from the magnitude of the denominators, this would only decrease the memory requirements by about 10% from that of going to the next HWM.

It is worth considering the appearance of, for example, the denominator of the Convergent Before the 2nd Generation HWM at Coefficient Number 9827. There is a 47 digit long Preamble, followed by 62,999,953 '9's, followed by a 46 digit long Penultimate, followed by 488,883 '0's. Files containing the 2nd Generation HWM Convergent denominators may be found in [16][17].

# 3rd, 4th and Higher Generation HWMs

As the HWMs get larger, more Generations of HWMs begin to appear between the HWMs. The table below was constructed by looking at the $C_{10}$ CFE coefficient lengths. Included in the table are those coefficients greater than 50 digits between HWM #5 and HWM #9, greater than 300 digits between HWM #9 and HWM #10, greater than 5000 digits between HWM #10 and HWM #11, and greater than 50000 digits between HWM #11 and HWM #12. The limits prevent HWM Generation Numbers above 4th Generation from appearing, but are not so restrictive as to hide other large coefficients.

### Table 6 - HWM Generations' Order and Lengths

| Coefficient Number | Number of Digits | | | |
|---|---|---|---|---|
| | 1st Generation HWM | 2nd Generation HWM | 3rd Generation HWM | 4th Generation HWM |
| 18 | 166 | | | |
| 40 | 2504 | | | |
| 101 | | 140 | | |
| 162 | 33,102 | | | |
| 246 | | | 109 | |
| 357 | | 2468 | | |
| 459 | | | 136 | |
| 526 | 411,100 | | | |
| 638 | | | | 90 |
| 820 | | | 2423 | |
| 1051 | | | | 63 |
| 1221 | | 33,056 | | |
| 1362 | | | | 95 |
| 1515 | | | 2458 | |
| 1627 | | | | 120 |
| 1708 | 4,911,098 | | | |
| 2074 | | | | 2411 |
| 2528 | | | 33,005 | |
| 3071 | | | | 2374 |
| 3569 | | 411,044 | | |
| 3916 | | | | 2419 |
| 4311 | | | 33,051 | |
| 4615 | | | | 2457 |
| 4838 | 57,111,096 | | | |
| 5810 | | | | 32,985 |
| 6992 | | | 410,979 | |
| 8469 | | | | 32,939 |
| 9827 | | 4,911,032 | | |
| 10,822 | | | | 32,996 |
| 12,015 | | | 411,036 | |
| 12,869 | | | | 33,041 |
| 13,522 | 651,111,094 | | | |

| 15,924 | | | | 410,961 |
| 18,884 | | | 4,910,961 | |
| 22,153 | | | | 410,907 |
| 25,069 | | 57,111,020 | | |
| 27,434 | | | | 410,973 |
| 30,231 | | | 4,911,027 | |
| 32,533 | | | | 411,024 |

The shaded rows are the rows with the HWMs. We feel that other than the 2nd Generation HWM lengths as covered in Conjecture 8, more data is needed to analyze the patterns of the lengths of the 3rd and 4th Generation HWMs. However, a pattern in the order in which the different Generations appear relative to each other is evident by looking at the table. It is also evident that there are higher levels of Generations of HWMs such as 5th Generation, etc.

# Section 4      Notes on Computing

This section gives more information on the various programs that were developed for, and used in the study. Information on the hardware that was used in the study is also given.

## Ruby Programs

As stated earlier, we initially used a Ruby program to calculate the HWM Convergents from the CFE coefficients available online up to Coefficient Number 161. The error of the Convergents was also calculated. This was relatively easy as arbitrary precision arithmetic functionality is built into the Ruby standard libraries.

We noticed patterns in the denominators of the HWM Convergents. A Ruby program was written to calculate, for each HWM Convergent, the Convergent's denominator using Conjecture 6. The program calculates the numerator by multiplying the denominator by a sufficient number of $C_{10}$ digits as determined from Conjecture 7. The program then calculates the CFE coefficients up to each HWM, and calculates the approximate error of the convergent.

The Convergent Before HWM #9 was the last one in which the error could practically be calculated with the Ruby program. The error calculation required approximately 12 times the number of digits as the CFE calculations and the construction of the $C_{10}$ Constant itself was time consuming for the Ruby program.

Therefore, the Ruby program was modified and the error was not calculated. This allowed for the calculation of the Convergent Before HWM #10. However, the calculation time was beginning to grow significantly, being approximately 140 times longer for HWM #10 than for HWM #9 (see Table 7). An increase in the amount of time that it took for the Ruby program to generate $C_{10}$ to a sufficient number of digits for calculating the numerator was noticed.

A program was written in C to construct and store to a file, the number of digits of $C_{10}$ required for the calculation of the numerator. The Ruby program was modified to read the C file instead of constructing $C_{10}$. This decreased the calculation time by about 15 % (much more so for the error calculating program) but the rate of increase in calculation time from one HWM to the next indicated that continuation to higher HWMs with the Ruby program was impractical.

## C Programs

We then switched to a C program not only to construct $C_{10}$ to the required number of digits, but also to perform all calculations. The MPIR 32 bit library (for use on a 64 bit machine) was built and linked in order for the C program to handle the required number of digits. The speed of the C program allowed for the calculation of the error for the Convergent Before HWM #10 and it also allowed for the calculation of the Convergent Before HWM #11. However, for the Convergent Before HWM #12, almost 790 million decimal digits are needed. This is currently beyond the capability of MPIR.

As a result, a change to the 64 bit GMP library (and a different OS) was made to calculate the Convergent before HWM #12. This calculation computed HWM #11, the largest which has been calculated to date by other means [15]. Furthermore, the Coefficient Number of HWM #12 (34062) was calculated, a result which is believed to be previously unknown. The GMP library was also used to verify the error for the Convergent Before HWM #11.

## Calculation Details

The following table shows the computation time of the various methods for calculating the Convergents Before HWM #N. All results were obtained on a desktop computer with a 3.40 GHz i7-2600 CPU, running single threaded, with 8 GB of 1333 MHz DDR3 SDRAM. The programs write the following to a file: the CFE coefficients as calculated from the Convergent, the numerator of the Convergent, and where applicable, the value of $C_{10}$ as calculated by the Convergent, and the approximate error from $C_{10}$.

The following table gives the calculation and storage run times:

**Table 7 – Calculation Times**

| Convergent Before HWM Number (N) | Time to Run, Including Writing to Files, Seconds | | | | | | | |
| --- | --- | --- | --- | --- | --- | --- | --- | --- |
| | With Calculating Error | | | | Without Calculating Error | | | |
| | Ruby Program, Self Generating $C_{10}$ | Ruby Program, $C_{10}$ Read From C Generated File | MPIR C Program | GMP C Program | Ruby Program, Self Generating C10 | Ruby Program, $C_{10}$ Read From C Generated File | MPIR C Program | GMP C Program |
| 7 | 0.22 | 0.11 | | | 0.03 | 0.03 | | |
| 8 | 9.15 | 1.88 | 1 | | 0.67 | 0.61 | < 1 | |
| 9 | 1145 | 113 | 62 | 9 | 75 | 67 | 3 | |
| 10 | | | 642 | 158 | 10,457 | 9610 | 62 | 16 |
| 11 | | | | 2769 | | | 1340 | 381 |
| 12 | | | | | | | | 10,692 |

Programs were written in Ruby and C to calculate the 2nd Generation HWM Convergents (from the CFE coefficients calculated using the Convergents of the HWMs) and the error for each such Convergent. The performance of these programs is not tabulated in this paper.

## Programs to Check the Results

Programs were also written to check the results of the various programs against each other to make sure that they matched, and to find the lengths of the coefficients, etc.

## Location of the Programs and Output Files

The programs may be found at http://code.google.com/p/champernowne-constant-cfe-coefficient-calculation-hwm/downloads/list. The programs are available under the GNU LGPL (GNU Lesser General Public License). The documentation consists of a readme.txt file and program comments. The resulting data files are present as well in the various subfolders (through the Convergent before HWM #11), so that researchers do not need to perform the calculations to see the results. The file containing the coefficients calculated from the Convergent Before HWM #12 may be found at https://docs.google.com/file/d/0B_ZIuCD9HzQtNnlLQUdIYTRXY2M/edit (it is about 789 MB).

# Section 5        Suggestions for Future Analysis

In this section we give suggestions for improvements and extensions of the present work.

## Computing

The programs that have been written for the current analysis have focused on improving the execution time (Ruby to C) and on overcoming an inherent limitation in the size of the numbers (MPIR to GMP) that may be processed. The limitation for the next calculation, the Convergent Before HWM #13, is memory space. The current program is not efficient from a memory standpoint. Therefore, the next step in the future research of this subject should be to optimize the program to be as memory efficient as possible.

Nevertheless, for researchers that have access to machines with greater memory, it is conceivable that the Convergent before HWM #13 is within reach. This computation would reveal the value of the predicted 7,311,111,092 digit long HWM #12, and it would also give the Coefficient Number of HWM #13.

It may also be possible to utilize multi-threading on different cores of the processor to increase the speed and reduce calculation time. Doing so would probably be of greatest advantage during the calculation of the numerator as the results from the various threads could be added together after processing to get the correct value for the numerator. The additions must be performed using floating point operations and the Ceil function should be taken on the result.

## HWM Coefficient Number Prediction or Pattern

A pattern in the HWM Coefficient Numbers has not been observed. It may be possible for future researchers to find such a pattern.

## $2^{nd}$ Generation HWM Convergent Denominators: Formulae for the Preamble and Penultimate

A pattern for the general construction of the denominators of the Convergents of the $2^{nd}$ Generation HWMs has been found. The total length of each denominator, the length of each Preamble, the number of nines following each Preamble, the length of each Penultimate, and the number of zeroes following each Penultimate is known. However, the patterns for the exact values of each Preamble and Penultimate have not been found. If such patterns are found, it would be possible to calculate the $C_{10}$ CFE Coefficients Before the $2^{nd}$ Generation HWMs. This would require approximately 10% less memory, and less computing time than computing the Coefficients Before the HWMs.

## Patterns in the HWM Convergent Characteristics of Champernowne's Constant in Bases Other Than 10

An investigation to determine if similar patterns and characteristics exist in the HWM Convergents of $C_X$, where $C_X$ is Champernowne's Constant in base X, could prove to be a worthwhile project. It is assumed that base 10 is not a special case. Therefore patterns and characteristics similar to those found in base 10 are likely to exist.

## Possible Proof of the Conjectures

It is hoped that some, if not all, of the conjectures will eventually be proven to be true.

# Conclusions

A very efficient method for calculating the CFE coefficients for Champernowne's Constant in Base 10 ($C_{10}$) has been presented and the current findings have been given. It is conjectured that this method calculates the coefficients up to the High Water Mark (HWM) Coefficients for HWM Numbers greater than or equal to 4. It was pointed out that the primary disadvantage is that each successive HWM calculation requires an increase of approximately 12 times in memory and an increase of approximately 24 times in computing time.

It was stated that the calculation of the Convergent Before HWM #12 verified the largest known HWM (HWM #11), and also determined the Coefficient Number of HWM #12, which is believed to be previously unknown.

In addition, a number of conjectures have been stated that point out patterns in the HWM lengths, and the HWM Convergent properties.

The other large CFE coefficients, that are not themselves HWMs, have been coined as 2[nd] Generation, 3[rd] Generation, etc., HWMs. Conjectures have been presented which predict the lengths of the 2[nd] Generation HWMs, and predict patterns in the structure of the denominators for the 2[nd] Generation HWM Convergents. A table was presented that shows the distribution of the 1[st] through 4[th] Generations of HWMs.

Finally, a Summary of Findings is given on the following page.

# Summary of Findings

The following is a bulleted summary of the findings, predictions, and conjectures of the paper. For the details and the applicable formulae, please refer to the previous sections.

- The Coefficient Number of HWM #12 has been calculated to be 34,062. The exact value of HWM #12 remains unknown, and it is predicted to have exactly 7,311,111,092 digits (Table 1).
- The Convergents of the $C_{10}$ CFE, truncated immediately before the HWMs, for HMW #N, N ≥ 4, calculate $C_{10}$ accurately up to the last 8 before an integer power of 10 in the sequence used to construct $C_{10}$ (Conjecture 1).
- An equation is given for the error calculated by the Convergents Before HWM #N (Conjecture 4).
- An equation is given for the denominators of the Convergents Before HWM #N (Conjecture 6).
- The number of digits of $C_{10}$ that are required to calculate each of the numerators of the Convergents Before HWM #N is given (Conjecture 7).
- The CFE coefficients may be efficiently calculated using the numerators and denominators of the Convergents Before HWM #N (Conjecture 7).
- An equation is given for the lengths in decimal digits of the HWM CFE coefficients of $C_{10}$ (Conjecture 5).
- The HWM Coefficient Numbers must be even per the coefficient numbering scheme given in this paper (Conjecture 3).
- The terms Child HWM, or 2nd, 3rd, etc. Generation HWM have been introduced.
- An equation is given for the lengths of the 2nd Generation HWMs (Conjecture 8).
- An equation is given for the error calculated by the Convergents Before the 2nd Generation HWMs (Conjecture 9).
- The 2nd Generation HWM coefficient numbers must be odd per the coefficient numbering scheme given in this paper (Conjecture 9).
- An equation is given for the length in decimal digits and the general structure is defined for the denominators of the Convergents Before the 2nd Generation HWMs (Conjecture 10).
- The distribution of the 2nd, 3rd, and 4th Generation HWMs has been tabulated up to HWM #12 (Table 6).
- A number of programs are available under the LGPL for calculating the findings presented in this paper. These programs may be downloaded from http://code.google.com/p/champernowne-constant-cfe-coefficient-calculation-hwm/downloads/list.
- The following are not discussed in the text of the paper:
  - Starting at HWM #5, the numerators of the Convergents before the odd numbered HWMs each has a somewhat large number of consecutive nines. The numerators for the Convergents before HWMs #5, #7, #9, and #11 have 5, 173, 2869, and 35987 consecutive nines, respectively. The numerator for the Convergent before HWM #11 also has a separate string of 165 consecutive nines.

- Starting with the Convergent Before HWM #6, the numerator for the Convergent Before HWM #N ends in '40…09', where the number of zeroes is (N - 5). This is used as one of the checks for a calculation.

*E-mail Address:*  john.sikora@xtera.com